\documentclass{birkjour}
\usepackage [latin1]{inputenc}

\newtheorem{thm}{Theorem}[section]
\newtheorem{cor}[thm]{Corollary}
\newtheorem{lem}[thm]{Lemma}

\theoremstyle{definition}

\theoremstyle{remark}
\newtheorem{rem}[thm]{Remark}

\numberwithin{equation}{section}

\newcommand{\BibTeX}{B\kern-0.1emi\kern-0.017emb\kern-0.15em\TeX}
\newcommand{\XYpic}{$\mathrm{X\kern-0.3em\raisebox{-0.18em}{Y}}$-$\mathrm{pic}\,$}

\newcommand{\cl}{C \kern -0.1em \ell}  

\newcommand{\ed}{\end{document}}

\begin{document}

\setcounter{page}{1} \setlength{\unitlength}{1mm}\baselineskip
.58cm \pagenumbering{arabic} \numberwithin{equation}{section}

\title[A note on gradient Solitons]
{A note on gradient Solitons on two classes of almost Kenmotsu Manifolds}

\author[ K. De and  U. C. De ]
{ Krishnendu De$^{*}$ and Uday Chand De }

\address
 {Department of Mathematics,
 Kabi Sukanta Mahavidyalaya,
The University of Burdwan.
Bhadreswar, P.O.-Angus, Hooghly,
Pin 712221, West Bengal, India. ORCID iD: https://orcid.org/0000-0001-6520-4520}
\email{krishnendu.de@outlook.in }

\address
{Department of Mathematics, University of Calcutta, West Bengal, India. ORCID iD: https://orcid.org/0000-0002-8990-4609}
\email {uc$_{-}$de@yahoo.com}
\footnotetext {$\bf{2020\ Mathematics\ Subject\ Classification\:}.$ 53D50, 53C25, 53C80.
\\ {Key words and phrases: $(k,\mu)'$- almost Kenmotsu manifolds, Kenmotsu manifolds, $(m,\rho)$-quasi Einstein solitons.\\
\thanks{$^{*}$ Corresponding author}
}}
\maketitle

\vspace{1cm}

\begin{abstract}
The purpose of the article is to characterize \textbf{gradient $(m,\rho)$-quasi Einstein solitons} within the framework of two classes of almost Kenmotsu Manifolds. Finally, we consider an example to justify a result of our paper.
\end{abstract}

\maketitle

	\section{{Introduction}}

The techniques of contact geometry carried out a significant role in contemporary mathematics and  consequently it have become well-known among the eminent researchers. Contact geometry has manifested from the mathematical formalism of classical mechanics. This subject matter has more than one attachments with the alternative regions of differential geometry and outstanding applications in applied fields such as phase space of dynamical systems, mechanics, optic and thermodynamics. In the existing paper, we have a look at the gradient solitons with nullity distribution which play a functional role in coeval mathematics.\par

In the study of Riemannian manifolds $(M,g)$, Gray \cite{ga} and Tanno \cite{ta} introduced the concept of \emph{$k$-nullity distribution $(k\in\mathbb{R})$} and is defined for any $p\in M$ and $k\in\mathbb R$ as follows:
\begin{equation}\label{a1}
N_{p}(k)=\{W\in T_{p}M:R(U,V)W = k[g(V,W)U-g(U,W)V]\},
\end{equation}
for any $U,V\in T_{p}M$, where $R$ indicates the Riemannian curvature tensor of type $(1,3)$.\par
Recently, the generalized idea of the $k$-nullity distribution, named as \emph{$(k,\mu)$-nullity distribution} on a contact metric manifold $(M^{2n+1},\eta, \xi, \phi,g)$ intro by  Blair, Koufogiorgos and Papantoniou \cite{bkp} and defined for any $p\in M^{2n+1}$ and $k,\mu\in\mathbb R$ as :
\begin{eqnarray}\label{a2}
N_{p}(k,\mu)=\{W\in T_{p}M^{2n+1}:R(U,V)W&=&k[g(V,W)U-g(U,W)V]\nonumber\\&&+\mu[g(V,W)h U-g(U,W)h V]\},
\end{eqnarray}
for any $U,V\in T_{p}M$ and $h=\frac{1}{2}\pounds_{\xi}\phi$, where $\pounds$ denotes the Lie differentiation.\par

In $2009$, Dileo and Pastore \cite{dp2} introduced a further generalized concept of the $(k,\mu)$-nullity distribution, called the \emph{$(k,\mu)'$-nullity distribution} on an almost Kenmotsu manifold $(M^{2n+1},\eta,\xi,\phi, g)$ and is defined for any $p\in M^{2n+1}$ and $k,\mu\in\mathbb R$ as :
\begin{eqnarray}\label{a3}
N_{p}(k,\mu)'=\{W\in T_{p}M^{2n+1}:R(U,V)W&=&k[g(V,W)U-g(U,W)V]\nonumber\\&&+\mu[g(V,W)h'U-g(U,W)h'V]\},
\end{eqnarray}
for any $U,V\in T_{p}M$ and $h'=h\circ\phi$.\par

A Riemannian metric $g$ on a Riemannian manifold $M$ is called a \emph{gradient $(m,\rho)$-quasi Einstein soliton} if there exists a smooth function $f:M^{n}\rightarrow \mathbb{R}$ and three real constants $\rho$, $\lambda$ and $m$ $(0<m\leq \infty)$ such that

\begin{equation} S+\nabla^{2}f-\frac{1}{m}df \otimes df =\beta g =(\rho r+\lambda) g,\label{a8}\end{equation}

where $\nabla^{2}$ and $\otimes$ indicate the Hessian of $g$ and tensor product, respectively. The expression $S+\nabla^{2}f-\frac{1}{m}df \otimes df$ is the $m$-Bakry-Emery Ricci tensor, which is proportional to the metric $g$ and $\lambda = constant$ \cite{ww}.
The soliton becomes trivial if the potential function $f$ is constant and the triviality condition implies that the manifold is an Einstein manifold. Furthermore, when $m=\infty$, the foregoing equation reduces to gradient $\rho$-Einstein soliton. This notion was introduced in \cite {cm} and recently, Venkatesha et al. studied $\rho$- Einstein solitons \cite{ven} on almost Kenmotsu manifold. In this connection, the properties of $(m,\rho)$-quasi Einstein solitons in different geometrical structures have been studied (in details) by (\cite{ag1}, \cite{hw}) and others. To know more about almost Kenmotsu manifold, here we may mention the work of Venkatesha et al. and Wang and his collaborator (\cite{ven1}, \cite{wang1}, \cite{wang2}, \cite{wa2}). Recently, Wang \cite{wang} has studied Yamabe solitons and gradient Yamabe solitons within the context of almost Kenmotsu $(k,\mu)'$ manifolds.\par
Motivated from the above studies, we make the contribution to investigate gradient $(m,\rho)$-quasi Einstein solitons in almost Kenmotsu Manifolds.\par

The current paper is constructed as :\
In section 2, we recall some basic facts and formulas of almost Kenmotsu manifolds which we will need throughout the paper. In sections 3, we characterize gradient $(m,\rho)$-quasi Einstein solitons with nullity distribution and obtain some interesting results. In the next section we consider the gradient $(m,\rho)$-quasi Einstein soliton in a $3$-dimensional Kenmotsu manifold and proved that the manifold is of constant sectional curvature $-1$, provided $r\neq-(\lambda+2)$. Then we consider an example to verify the result of our paper.

\section{{Almost Kenmotsu manifolds}}

In this section we gather the formulas and results of almost Kenmotsu manifolds which will be required on later sections.
A differentiable manifold $M^{2n+1}$ of dimension $(2n+1)$ is called \emph {almost contact metric manifold} if it admits  a covariant vector field $\eta$, a contravariant vector field $\xi$, a $(1,1)$ tensor field $\phi$ and a Riemannian metric $g$ such that
\begin{equation}\label{b1}
\phi^{2}=-I+\eta\otimes\xi,\;\eta(\xi)=1,
\end{equation}
\begin{equation}\nonumber
g(\phi U,\phi V)=g(U,V)-\eta(U)\eta(V),
\end{equation}
where $I$ indicates the identity endomorphism . Then also $\phi\xi=0$ and $\eta\circ\phi=0$; in a straight forward calculation both can be extracted from (\ref{b1}).\par
On an almost Kenmotsu manifold $M^{2n+1}$, the two symmetric tensor fields $h=\frac{1}{2}\pounds_{\xi}\phi$ and $l=R(\cdot,\xi)\xi$, satisfy the following relations \cite{dp2}
\begin{equation}\label{b2}
h\xi=0,\;l\xi=0,\;tr(h)=0,\;tr(h')=0,\;h\phi+\phi h=0,
\end{equation}
\begin{equation}\label{b3}
\nabla_{U}\xi=-\phi^2U+h'U(\Rightarrow \nabla_{\xi}\xi=0),
\end{equation}
\begin{equation}\label{b4}
\phi l\phi-l=2(h^{2}-\phi^{2}),
\end{equation}
\begin{equation}\label{b5}
R(U,V)\xi=\eta(U)(V-\phi hV)-\eta(V)(U-\phi hU)+(\nabla_{V}\phi h)U-(\nabla_{U}\phi h)V,
\end{equation}
for any vector fields $U,V$.\par
Let the Reeb vector field $\xi$ of an almost Kenmotsu manifold belonging to the $(k,\mu)'$-nullity distribution.  Then the symmetric tensor field $h'$ of type $(1,1)$ satisfies the relations $h'\phi+\phi h'=0$ and $h'\xi=0$.
 Also, it is understandable that
\begin{equation}\label{b6}
h'^{2}=(k+1)\phi^2\, (\Leftrightarrow h^{2}=(k+1)\phi^2),\;\;h=0\Leftrightarrow h'=0.
\end{equation}
From the equations (\ref{b1}) and (\ref{b6}), it follows that $k\leq-1$.
Certainly, by (\ref{b6}), we
 conclude that  $h'$ vanishes if and only if $k=-1$.
In view of the  Proposition 4.1 of \cite{dp2} on a $(k,\mu)'$-almost Kenmotsu manifold with $k <-1$ leads to $\mu=-2$.

For an almost Kenmotsu manifold, we possess from (\ref{b3})
\begin{equation}\label{b7}
R(U,V)\xi=k[\eta(V)U-\eta(U)V]+\mu[\eta(V)h'U-\eta(U)h'V],
\end{equation}
\begin{equation}\label{b8}
R(\xi,U)V=k[g(U,V)\xi-\eta(V)U]+\mu[g(h'U,V)\xi-\eta(V)h'U],
\end{equation}
where $k,\mu\in\mathbb R$. Contracting $V$ in (\ref{b8}) we have
\begin{equation}\label{b9}
S(U,\xi)=2k\eta(U).
\end{equation}
Also, the distribution which is indicated by $\mathcal{D}$ (defined by $\mathcal{D}=$ker $\eta$). Presume $X\in \mathcal{D}$ be the eigen vector of $h'$ corresponding to the eigen value $\lambda$. Then $\lambda^{2}=-(k+1)$, a constant, which follows from (\ref{b6}). Therefore $k\leq -1$ and $\lambda=\pm\sqrt{-k-1}$. The non-zero eigen value $\lambda$ and $-\lambda$ are respectively indicated by $[\lambda]'$ and $[-\lambda]'$, which are the corresponding eigen spaces associated with $h'$. Now before introducing the detailed proof of our main theorem, we first write the following lemma:
\begin{lem} \label{L2}(Lemma. 3.2 of \cite{wal}) Let $(M^{2n+1},\eta,\xi,\phi,g)$ be an almost Kenmotsu manifold such that $\xi$ belongs to the $(k,\mu)'$-nullity distribution and $h'\neq0$.
Then the Ricci operator $Q$ of $M^{2n+1}$ is given by
\begin{eqnarray} \label{b10} Q =-2nid+2n(k+1)\eta \otimes \xi - 2nh',\end{eqnarray}
where $k < -1$, moreover, the scalar curvature of $M^{2n+1}$ is $2n(k-2n)$.

\end{lem}

\section{Gradient $(m,\rho)$-quasi Einstein solitons on a $(2n+1)$-dimensional almost Kenmotsu manifold}

Let us assume that the Riemannian metric of a $(2n+1)$-dimensional almost Kenmotsu manifold such that $\xi$ belongs to the $(k,\mu)'$-nullity distribution and $h'\neq0$, is a $(m,\rho)$-quasi Einstein soliton. Then the equation (\ref{a8}) may be expressed as

\begin{equation}\label{c31}
\nabla_{U}D f+Q U=\frac{1}{m}g(U,D f)D f+\beta U.
\end{equation}

Taking covariant derivative of (\ref{c31}) along the vector field $V$, we get
\begin{eqnarray}\label{c32}
\nabla_{V}\nabla_{U}D f&=&-\nabla_{V}QU+ \frac{1}{m}\nabla_{V}g(U,D f)D f\nonumber\\&&+\frac{1}{m}g(U,D f)\nabla_{V}D f+\beta \nabla_{V}U.
\end{eqnarray}
Interchanging $U$ and $V$ in (\ref{c32}), we lead
\begin{eqnarray}\label{c33}
\nabla_{U}\nabla_{V}D f&=&-\nabla_{U}QV+ \frac{1}{m}\nabla_{U}g(V,D f)D f\nonumber\\&&+\frac{1}{m}g(V,D f)\nabla_{U}D f+\beta \nabla_{U}V
\end{eqnarray}
and
\begin{eqnarray}\label{c34}
\nabla_{[U,V]}D f=-Q[U,V]+ \frac{1}{m}g([U,V],D f)D f+\beta [U,V].
\end{eqnarray}
Equations (\ref{c31})-(\ref{c34}) and the symmetric property of Levi-Civita connection together
with $R(U,V)D f=\nabla_{U}\nabla_{V}D f-\nabla_{V}\nabla_{U}D f-\nabla_{[U,V]}D f$ we lead
\begin{eqnarray}\label{c35}
R(U,V)D f&=& (\nabla_{V}Q)U-(\nabla_{U}Q)V+\frac{\beta}{m}\{ (V f)U-(U f)V\}\nonumber\\&&
+\frac{1}{m}\{ (U f)QV-(V f)QU \}.
\end{eqnarray}
Taking inner product of (\ref{c35}) with $\xi$, we have
\begin{eqnarray}\label{c39}
g(R(U,V)D f,\xi)&=&\frac{\beta}{m}\{ (V f)\eta (U)-(U f)\eta(V)\}\nonumber\\&&
+\frac{1}{m}\{ (U f)\eta(QV)-(V f)\eta(QU) \}.
\end{eqnarray}
Again (\ref{b7}) implies that
\begin{equation}\label{c8}
g(R(U,V)\xi,Df)=k[\eta(V)(U f)-\eta(U)(V f)]+\mu[\eta(V)(h'U f)-\eta(U)(h'V f)].
\end{equation}
Combining equation (\ref{c39}) and (\ref{c8}) reveal that
\begin{eqnarray}\label{c40}
&&k[\eta(U)(V f)-\eta(V)(U f)]+\mu[\eta(U)(h'V f)-\eta(V)(h'U f)]\nonumber\\&&
=\frac{\beta}{m}\{ (V f)\eta (U)-(U f)\eta(V)\}\nonumber\\&&
+\frac{1}{m}\{ (U f)\eta(QV)-(V f)\eta(QU) \}.
\end{eqnarray}
Setting $V=\xi$ in the foregoing equation yields
\begin{equation} \label{c41} 2mh'D f=(\beta-2nk-mk)\{D f-(\xi f)\xi\},\end{equation}
where we have used $\mu=-2$.\par
Letting $\beta=(2nk+mk+2m)$ and taking into account the equation (\ref{b6}) and operating $h'$ on (\ref{c41}) produces that
\begin{equation} \label{c42} -(k+1)\{D f-(\xi f)\xi\}=h'D f,\end{equation}
Comparing the antecedent relation with (\ref{c41}) gives that
\begin{equation}\label{c43}
(k+2)\{D f-(\xi f)\xi\}=0.
\end{equation}
This shows that either $k=-2$ or $\{D f-(\xi f)\xi\}=0$.
\par
Case i: If $k=-2$, then the Proposition 4.1  and Corollary 4.2 of \cite{dp2} state that $M^{2n+1}$ is locally isometric to the Riemannian product $\mathbb{H}^{n+1}(-4)\times \mathbb{R}^{n}$.
 In fact, from (\cite{pet},\cite{pet1}) we can say that the product $\mathbb{H}^{n+1}(-4)\times \mathbb{R}^{n}$ is a rigid gradient Ricci soliton.\par
Case ii: \begin{equation} Df=(\xi f)\xi.\label{d11}\end{equation}
 Executing the covariant differentiation of (\ref{d11}) along $U\in \chi (M)$ and utilizing (\ref{b3}) we get
\begin{equation} \nabla _{U}Df=U(\xi f)\xi+(\xi f)U-(\xi f)\eta(U)\xi+(\xi f)h'U.\label{d12}\end{equation}
Replacing the foregoing equation into (\ref{c31}) revels that
\begin{equation} Q X=\beta U -(\xi f)U-U(\xi f)\xi+(\xi f)\eta(U)\xi-(\xi f)h'U+\frac{1}{m}g(U,D f)D f.\label{d13}\end{equation}
Comparing (\ref{b10}) and (\ref{d13}) give that
\begin{eqnarray} && \{\beta+2n -(\xi f)\}U+\{2n-(\xi f)\}h'U +\{(\xi f)\eta(U)\nonumber\\&&-U(\xi f)-2n(k+1)\eta(U)\}\xi+\frac{1}{m}g(U,D f)D f =0.\label{d14}\end{eqnarray}
Now operating $h'$ we get
\begin{eqnarray} && \{\beta+2n -(\xi f)\}h'U+(k+1)\{2n-(\xi f)\}\phi^{2}U \nonumber\\&&+\frac{1}{m}g(U,D f)h'D f =0.\label{d15}\end{eqnarray}
Contracting $U$ in the above equation we infer
\begin{equation} 2n(k+1)\{2n-(\xi f)\}+\frac{1}{m}g(h'D f,D f)=0.\label{d16}\end{equation}
Putting $U=\xi$ in (\ref{d15}), we have
\begin{equation} \eta(D f)h' D f=0.\label{d17}\end{equation}
This shows that either $\eta(D f)=0$ or $h' D f=0.$\par
Case (i): If $\eta(D f)=0$, then we get $\xi f=0$ and hence from equation (\ref{d11}) it follows that $f=constant$. Then we get from (\ref{c31}) that the manifold is an Einstein manifold.\par
Case (ii): If $h' D f=0$, then from (\ref{d16}) we get $(\xi f)=2n$. Now putting this value in (\ref{d15}) we obtain $\beta h' U=0$. Hence either $\beta=0$, or $h' U=0$. Both contradicts our assumptions $\beta=(2n+m)k+2m$ and $k<-1$  respectively .\par
Thus, we can state the following theorem:
\begin{thm}\label{main3}
Let the Riemannian metric of a $(2n+1)$-dimensional $(k,\mu)'$ almost Kenmotsu manifold with $h'\neq0$, be the gradient $(m,\rho)$-quasi Einstein soliton. Then either $M^{2n+1}$ is locally isometric to a rigid gradient Ricci soliton $\mathbb{H}^{n+1}(-4)\times \mathbb{R}^{n}$ or $M^{2n+1}$ is an Einstein manifold, provided $\beta=(2nk+mk+2m)$.
\end{thm}
We know that when $m=\infty$, the $(m,\rho)$-quasi Einstein soliton becomes gradient $\rho$-Einstein soliton. Putting the value $m=\infty$ in (\ref{c40}) and by a straight forward calculation we find that the manifold $M^{2n+1}$ is locally isometric to a rigid gradient Ricci soliton $\mathbb{H}^{n+1}(-4)\times \mathbb{R}^{n}$. Thus, we can state:

\begin{cor}\label{cor1}
Let the Riemannian metric of a $(2n+1)$-dimensional almost Kenmotsu manifold such that $\xi$ belongs to the $(k,\mu)'$-nullity distribution and $h'\neq0$, be the gradient $\rho$- Einstein soliton. Then $M^{2n+1}$ is locally isometric to a rigid gradient Ricci soliton $\mathbb{H}^{n+1}(-4)\times \mathbb{R}^{n}$.
\end{cor}
\begin{rem}
The above corollary have been proved by Venkatesha and Kumara in their paper \cite{ven}. They also prove that the potential vector field is tangential to the Euclidean factor $\emph{R}^{n}$.
\end{rem}

\section{Gradient $(m,\rho)$-quasi Einstein solitons on a $3$-dimensional Kenmotsu manifold}
In \cite{dp2}, Dileo and Pastore proved that the following conditions are equivalent:\par
(1) the (1,1)-type tensor field $h$ vanishes and foliations of the distribution $\mathcal{D}$ are K$\ddot{a}$hlerian.\par
(2) almost contact metric structure of an almost Kenmotsu manifold is normal.\par As a outcome, we get instantly that a 3-dimensional almost Kenmotsu manifold reduces to a Kenmotsu manifold if and only if $h=0$. In this section, we target to investigate the gradient $(m,\rho)$-quasi Einstein soliton on a $3$-dimensional Kenmotsu manifold. Making use of $h=0$ in equation (\ref{b3})
we get $\nabla \xi=-\phi^{2}$, and this revels that
\begin{equation} \label{k1} R(U,V)\xi =-\eta (V)U + \eta (U)V\end{equation}
for any $U,V$ $\in \chi(M)$ and therefore by contracting $V$ in (\ref{k1}) we obtain $Q\xi=-2\xi$.
From \cite{de} we know that for a 3-dimensional Kenmotsu manifold
\begin{eqnarray}\label{kk1} R(U,V)W &=&(\frac{r+4}{2})[g(V,W)U-g(U,W)V]\\&&
-(\frac{r+6}{2})[g(V,W)\eta (U)\xi-g(U,W)\eta (V)\xi\nonumber\\&&+\eta (V)\eta (W)U-\eta (U)\eta
(W)V]\nonumber,\end{eqnarray}
\begin{equation} \label{kk2} S(U,V)=\frac{1}{2}[(r+2)g(U,V)-(r+6)\eta
(U)\eta (V)].\end{equation}

Now before introducing the detailed proof of our main theorem, we first state the following result:
\begin{lem}\label{L3}(lemma. 4.1 of \cite{wa})
For a 3-dimensional Kenmotsu manifold $(M^{3},\phi,\xi,\eta,g)$, we have
\begin{eqnarray} \label{b13}
\xi r=-2(r+6)
\end{eqnarray}
\end{lem}
where $r$ denotes the scalar curvature of $M$.\\

Let us assume that the Riemannian metric of a $3$-dimensional almost Kenmotsu manifold be a $(m,\rho)$-quasi Einstein soliton. Then the equation (\ref{a8}) may be expressed as

\begin{equation}\label{k2}
\nabla_{U}D f+Q U=\frac{1}{m}g(U,D f)D f+\beta U.
\end{equation}
Taking covariant derivative of (\ref{k2}) along the vector field $V$, we get
\begin{eqnarray}\label{k3}
\nabla_{V}\nabla_{U}D f&=&-\nabla_{V}QU+ \frac{1}{m}\nabla_{V}g(U,D f)D f\nonumber\\&&+\frac{1}{m}g(U,D f)\nabla_{V}D f+\beta \nabla_{V}U+(V\beta)U.
\end{eqnarray}
Interchanging $U$ and $V$ in (\ref{k3}), we lead
\begin{eqnarray}\label{k4}
\nabla_{U}\nabla_{V}D f&=&-\nabla_{U}QV+ \frac{1}{m}\nabla_{U}g(V,D f)D f\nonumber\\&&+\frac{1}{m}g(V,D f)\nabla_{U}D f+\beta \nabla_{U}V+(U\beta)V
\end{eqnarray}
and
\begin{eqnarray}\label{k5}
\nabla_{[U,V]}D f=-Q[U,V]+ \frac{1}{m}g([U,V],D f)D f+\beta [U,V].
\end{eqnarray}
Equations (\ref{k2})-(\ref{k5}) and the symmetric property of Levi-Civita connection together
with $R(U,V)D f=\nabla_{U}\nabla_{V}D f-\nabla_{V}\nabla_{U}D f-\nabla_{[U,V]}D f$ we lead
\begin{eqnarray}\label{k6}
R(U,V)D f&=& (\nabla_{V}Q)U-(\nabla_{U}Q)V+\frac{\beta}{m}\{ (V f)U-(U f)V\}\nonumber\\&&
+\frac{1}{m}\{ (U f)QV-(V f)QU \}+\{(U\beta)V-(V\beta)U\}.
\end{eqnarray}
Taking inner product of (\ref{k6}) with $\xi$ and using (\ref{kk2}) we have
\begin{eqnarray}\label{k7}
g(R(U,V)D f,\xi)&=&\frac{\beta}{m}\{ (V f)\eta (U)-(U f)\eta(V)\}\nonumber\\&&
+\frac{1}{m}\{ (U f)\eta(QV)-(V f)\eta(QU) \}\nonumber\\&&
+\{(U\beta)\eta (V)-(V\beta)\eta(U)\}.
\end{eqnarray}
Again, from (\ref{k1}) we infer that
\begin{equation}\label{k8}
g(R(U,V)D f,\xi)=-\{ (V f)\eta (U)-(U f)\eta(V)\}.\end{equation}
Combining equation (\ref{k7}) and (\ref{k8}) reveal that
\begin{eqnarray}\label{k9}
-\{ (V f)\eta (U)-(U f)\eta(V)\}&=&\frac{\beta}{m}\{ (V f)\eta (U)-(U f)\eta(V)\}\nonumber\\&&
+\frac{1}{m}\{ (U f)\eta(QV)-(V f)\eta(QU) \}\nonumber\\&&+\{(U\beta)\eta (V)-(V\beta)\eta(U)\}.
\end{eqnarray}
Replacing $V$ by $\xi$ in the foregoing equation, we get
\begin{equation}\label{k10} d(f-\beta)=\xi(f-\beta)\eta,\end{equation}
where $d$ stands for the exterior differentiation, provided $\beta=-2$. In other word, $f-\beta$ is invariant along $\mathcal{D}$, i.e., $X(f-\beta)=0$ for any $X \in \mathcal{D}$. Taking into account the above fact and using Lemma 6.1, we infer
\begin{equation}\label{k11} (\xi f)=(\xi \beta)=\rho (\xi r)=-2\rho (r+6).\end{equation}
The contraction of the equation (\ref{k6}) along $U$ and applying Lemma 6.1, we get
\begin{eqnarray}\label{k12}
S(V,D f)&=& \frac{1}{2}(V r)+\frac{2 \beta }{m}(V f)\nonumber\\&&
-\frac{1}{m}\{r(V f)-g(Q V, D f)\}-2(V \beta).
\end{eqnarray}
Clearly, comparing the above equation with (\ref{kk2}) yields
\begin{eqnarray}\label{k13}
&& -(V r)-\frac{4 \beta }{m}(V f)
+\frac{2}{m}\{r(V f)-g(Q V, D f)\}\nonumber\\&&+4(V \beta)+(r+2)(V f)-(r+6)\eta (V)\xi f=0.
\end{eqnarray}
By a straight forward calculation, replacing $V$ by $\xi$ in (\ref{k13}) and using (\ref{k11}), we can easily obtain

\begin{equation} \label{k14} (\xi f)=\frac{2m(r+6)}{4\beta-2r-4}.\end{equation}
Comparing the antecedent equation with (\ref{k11}) reveals that
\begin{equation} \label{k15} 2(r+6)\{\rho+\frac{m}{4\beta-2r-4}\}=0.\end{equation}
This shows that either $r=-6$ or $\{\rho+\frac{m}{4\beta-2r-4}\}=0$.
Next, we split our investigation as :\par
Case (i): If $r=-6$, then from equation (\ref{kk2}) we conclude that $S=-2g$. Therefore, by using equation (\ref{kk1}) we sum up that
the manifold is of constant sectional curvature $-1$.\par
Case (ii): If $\{\rho+\frac{m}{4\beta-2r-4}\}=0$, then by a simple calculation we get
\begin{equation} \label{k16} r=\frac{4\lambda-4-4\rho\lambda-4\rho+m}{2(2\rho^{2}-3\rho+1)}=constant,\end{equation} provided $\rho\neq 1$. Hence by applying Lemma 6.1 we can easily get $r=-6$. Therefore, from Case (i) we see that the manifold is of constant sectional curvature $-1$. After combining the two conditions namely, $\beta=-2$ and $\rho\neq 1$, we can write $r\neq -(\lambda+2)$. Thus we have the following theorem:
\begin{thm}\label{main4}
Let the Riemannian metric of a $3$-dimensional Kenmotsu manifold be the gradient $(m,\rho)$-quasi Einstein soliton. Then the manifold is of constant sectional curvature $-1$, provided $r\neq -(\lambda+2)$.
\end{thm}
We know that when $m=\infty$, the $(m,\rho)$-quasi Einstein soliton gives the so called gradient $\rho$-Einstein soliton. Putting the value $m=\infty$ in (\ref{k9}) and by a straight forward calculation we find that the manifold is of constant sectional curvature $-1$. Thus, we can state:
\begin{cor}\label{cor2}
Let the Riemannian metric of a $3$-dimensional Kenmotsu manifold be the gradient $\rho$- Einstein soliton. Then the manifold is of constant sectional curvature $-1$.
\end{cor}

\section{\textbf{Example }}
Here we consider an example cited in our paper \cite{dek1}.

We consider the 3-dimensional manifold $M=\{(u,v,w)\varepsilon \mathbb{R}^{3}, w\neq0\},$ where
$(u,v,w)$ are standard coordinate of $\mathbb{R}^{3}.$

The vector fields$$E_{1}=w\frac{\partial }{\partial
u},\hspace{7pt}E_{2}=w\frac{\partial } {\partial v}
,\hspace{7pt}E_{3}=-w\frac{\partial }{\partial w}$$ are linearly
independent at each point of $M.$

Let $g$ be the Riemannian metric defined by
$$g(E_{1},E_{3})=g(E_{1},E_{2})=g(E_{2},E_{3})=0,$$
$$g(E_{1},E_{1})=g(E_{2},E_{2})=g(E_{3},E_{3})=1.$$

Let $\eta $ be the 1-form defined by $\eta (W)=g(W,E_{3})$ for any
$W\varepsilon \chi (M).$

Let $\phi $ be the $(1,1)$ tensor field defined by $$\phi
(E_{1})=-E_{2},\hspace{7pt} \phi (E_{2})=E_{1},\hspace{7pt}\phi
(E_{3})=0.$$
Then for $E_{3}=\xi $ , the structure $(\phi ,\xi ,\eta ,g)$
defines an almost contact metric structure on $M$.
Further Koszul's formula yields
$$\nabla _{E_{1}}E_{3}=E_{1},\hspace{10pt}\nabla _{E_{1}}E_{2}=0,\hspace{10pt}
\nabla _{E_{1}}E_{1}=-E_{3},$$
$$\nabla _{E_{2}}E_{3}=E_{2},\hspace{10pt}\nabla _{E_{2}}E_{2}=E_{3},\hspace{10pt}
\nabla _{E_{2}}E_{1}=0,$$ \begin{equation}\nabla
_{E_{3}}E_{3}=0,\hspace{10pt}\nabla _{E_{3}}E_{2}=0,\hspace{10pt}
\nabla _{E_{3}}E_{1}=0.\label{f5}\end{equation} From the above it follows
that the manifold satisfies $\nabla _{U}\xi=U-\eta (U)\xi$, for
$\xi=E_{3}$. Hence the manifold is a Kenmotsu manifold.\\
We verified that
$$R(E_{1},E_{2})E_{3}=0,\hspace{10pt}R(E_{2},E_{3})E_{3}=-E_{2},\hspace{10pt}
R(E_{1},E_{3})E_{3}=-E_{1},$$
$$R(E_{1},E_{2})E_{2}=-E_{1},\hspace{10pt}R(E_{2},E_{3})E_{2}=E_{3},\hspace{10pt}
R(E_{1},E_{3})E_{2}=0,$$
$$R(E_{1},E_{2})E_{1}=E_{2},\hspace{10pt}R(E_{2},E_{3})E_{1}=0,\hspace{10pt}
R(E_{1},E_{3})E_{1}=E_{3}.$$\\

From the above expressions of the curvature tensor $R$ we obtain
\begin{eqnarray}
S(E_{1},E_{1})&=&g(R(E_{1},E_{2})E_{2},E_{1})+g(R(E_{1},E_{3})E_{3},E_{1})
\nonumber\\&=&-2.\end{eqnarray}
Similarly, we have $$S(E_{2},E_{2})=S(E_{3},E_{3})=-2.$$
Therefore, $$r=S(E_{1},E_{1})+S(E_{2},E_{2})+S(E_{3},E_{3})=-6.$$
Let $f:M^{3}\rightarrow \mathbb{R}$ be a smooth function defined by $f=-w^{2}$.
Then gradient of $f$ with respect to $g$ is given by
$$  Df=-2w\frac{\partial }{\partial w}=2E_{3}.$$
With the help of (\ref{f5}) we can easily get
$$ Hessf(E_{3},E_{3})=0.$$
Thus gradient $\rho$-Einstein soliton equation shows
$$Hessf(E_{3},E_{3})+S(E_{3},E_{3})+2g(E_{3},E_{3})=0. $$
Similarly checking the other components we conclude that that $M^{3}$ satisfies
$$Hessf(U,V)+S(U,V)+2g(U,V)=0. $$
Thus $g$ is a gradient $\rho$-Einstein soliton with $f=-w^{2}$ and $\beta=-2$.
Hence the \textbf{Corollary 4.3} is verified.

\section{Declarations}
\subsection{Funding }
Not applicable.
\subsection{Conflicts of interest/Competing interests}
The authors declare that they have no conflict of interest.
\subsection{Availability of data and material }
Not applicable.
\subsection{Code availability}
Not applicable.

\end{document}